\documentclass[12pt,reqno]{amsart}   
\usepackage{amsthm,amsmath,amssymb,amsfonts,amscd,xypic,verbatim,latexsym}
\usepackage[mathscr]{eucal}
\newtheorem{Theorem}{Theorem}[section]
\newtheorem{Lemma}[Theorem]{Lemma}
\newtheorem{Proposition}[Theorem]{Proposition}
\newtheorem{Corollary}[Theorem]{Corollary}
\newtheorem{Remark}[Theorem]{Remark}

\newtheorem{Definition}[Theorem]{Definition}
\newtheorem{Example}[Theorem]{Example}

\renewcommand{\P}{{\mathbf P}} 
\renewcommand{\O}{{\mathcal O}} 
\newcommand{\B}{{\mathcal B}} 
\newcommand{\Le}{\mathcal L}  
\newcommand{\E}{{\mathcal E}} 
\newcommand{\Sc}{{\mathbb S}} 
\newcommand{\integers}{{\mathbb Z}}
\newcommand{\naturals}{{\mathbb N}}
\newcommand{\uL}{{\underline L}}
\newcommand{\ualpha}{{\underline \alpha}}
\newcommand{\uh}{{\underline h}}
\newcommand{\ux}{{\underline x}}
\newcommand{\uy}{{\underline y}}
\renewcommand{\phi}{\varphi}
\newcommand{\Sym}{\text{Sym}}
\newcommand{\image}{\text{image}\,}

\newcommand{\vspan}{\text{span}}
\renewcommand{\ker}{\text{ker}}
\newcommand{\rank}{\text{rank}}
\newcommand{\ann}{\text{ann}}

\newcommand{\Hom}{\text{Hom}}

\providecommand{\qed}{{\hfill $\Box$}}
\renewcommand{\qed}{{\hfill $\Box$}}
\newcommand{\demo}{{\sc proof.}\,\,\,}

\newcommand{\ra}{\rightarrow}
\newcommand{\RA}{\Rightarrow}

\newcommand{\lra}{\longrightarrow}

\newcommand{\lla}{\longleftarrow}

\begin{document}   
\title{On parameter spaces for Artin level algebras}   
\author{J.~V.~Chipalkatti and A.~V.~Geramita}
\maketitle   

\bigskip  
\parbox{12cm}{\small   
We describe the tangent space to the parameter variety of all artin level   
quotients of a polynomial ring in $n$ variables having specified socle  
degree and type. When $n=2$, we relate this variety  
to the family of secants of the rational normal curve. With additional  
numerical hypotheses, we prove a projective normality theorem for the  
parameter variety in its natural Pl{\"u}cker embedding.  }  
\bigskip  

Let $R=k[x_1,\dots,x_n]$ denote a polynomial ring and let   
$\uh: \naturals \lra  \naturals$ be a numerical function. Consider the  
set of all graded artin level quotients $A = R/I$ having Hilbert  
function $\uh$. This set (if nonempty) is naturally in bijection   
with the closed points of a quasiprojective scheme $\Le^\circ(\uh)$.   
The object of this note is to prove some specific geometric properties of these   
schemes, especially for $n=2$. The case of Gorenstein Hilbert functions   
(i.e., where $A$ has type $1$) has been extensively studied, and several   
qualitative and quantitative results are known (see \cite{IK}).   
Our results should be seen as generalizing some of them to the   
non-Gorenstein case.   

After establishing notation, we summarize the results in the next section.   
See \cite{Ger1,IK} as general references for most of the constructions   
used here.   
\section{Notation and Preliminaries}   
The base field $k$ will be algebraically closed of characteristic  
zero (but see Remark \ref{remark.char}).  
Let $V$ be an $n$-dimensional $k$-vector space, and let   
\[ R=\bigoplus_{i \ge 0} \Sym^i \, V^*,   
\quad S =\bigoplus_{i \ge 0} \Sym^i \, V. \]   
Let $\{x_1,\dots,x_n\},\{y_1,\dots,y_n\}$ be   
dual bases of $V^*$ and $V$ respectively, leading to identifications   
$R=k[x_1,\dots,x_n], S=k[y_1,\dots,y_n]$.   
There are internal products (see \cite[p.~476]{FH})
\[ \Sym^j\, V^* \otimes \Sym^i\, V \lra \Sym^{i-j}\, V, \quad   
u \otimes F \lra u.F \]   
making $S$ into a graded $R$-module. This action may be seen as  
partial differentiation; if $u(\ux) \in R$ and   
$F(\uy) \in S$, then    
\[ u.F = u(\partial/\partial y_1, \dots, \partial/ \partial y_n) \, F.   
\]
If $I \subseteq R$ is a homogeneous ideal, then $I^{-1}$ is the   
$R$-submodule of $S$ defined as $\{F \in S: u.F=0 \;\;   
\text{for all $u \in I$} \}$. This module (called Macaulay's inverse  
system for $I$) inherits a grading from $S$, thus   
$I^{-1} = \bigoplus\limits_i (I^{-1})_i$.  Reciprocally, if $M \subseteq S$   
is a graded submodule, then   
$\ann (M) = \{ u: u.F =0 \; \; \text{for all $F \in M$} \}$ is a homogeneous   
ideal in $R$.   
In classical terminology, if $u.F = 0$ and   
$\deg u \le \deg F$, then $u, F$ are said to be \emph{apolar} to each other.   

For any $i$, we have the Hilbert function   
$H(R/I,i) = \dim_k (R/I)_i = \dim_k (I^{-1})_i$.   
The following theorem is fundamental.   
\begin{Theorem}[Macaulay--Matlis duality]   
We have a bijective correspondence   
\[ \begin{aligned} \{ \text{homogeneous ideals $I < R$} \} &   
   \rightleftharpoons    \{ \text{graded $R$-submodules of $S$}\}  \\   
   I & \lra I^{-1} \\ \ann(M) & \lla  M \end{aligned}   
\]   
Moreover, $I^{-1}$ is a finitely generated $R$-module iff $R/I$ is artin.   

\end{Theorem}
Let $R/I=A$ be artin with graded decomposition   
\[ A= k \oplus A_1 \oplus \dots \oplus A_d, \quad A_d \neq 0, A_i=0 \;    
\text{for $i > d$}.  \]   
Recall that  
\[ \text{socle}(A) = \{ u \in A: u\, x_i =0 \; \; \text{for every $i$} \}.  
\]  
Then $A_d \subseteq \text{socle}(A)$, and $A$ is said to be \emph{level} if equality   
holds. This is true iff $I^{-1}$ is generated as an $R$-module by   
exactly $t:=\dim A_d$ elements in $S_d$. When $A$ is level, the number $t$ is   
called the \emph{type} of $A$, and it coincides with its Cohen-Macaulay type. Thus $A$   
is Gorenstein iff $t=1$. The number $d$ is the \emph{socle degree} of $A$.  
Altogether we have a bijection   
\[ \begin{aligned}   
\{A: \, & A=R/I \, \text{artin level of type $t$ and socle degree $d$} \}   
 \rightleftharpoons G(t,S_d);  \\    
A & \lra (I^{-1})_d, \quad \quad   
R/\ann(\Lambda)  \lla \Lambda.     
\end{aligned} \]   
Here $G(t,S_d)$ denotes the Grassmannian of $t$-dimensional vector  
subspaces of $S_d$. Notice the canonical isomorphism   
\begin{equation}  G(t,S_d) \simeq G(\dim R_d -t,R_d)   
\label{dual.grass} \end{equation}   
taking $\Lambda$ to $I_d$. We sometimes write $\Lambda_i$ for $(I^{-1})_i$ and   
$\Lambda$ for $\Lambda_d = (I^{-1})_d$.   
\begin{Remark} \rm   
The algebra $A$ is level iff $\Lambda_d$ generates $I^{-1}$ as an $R$-module,   
that is to say iff the internal product map   
\[ \alpha_i: R_{d-i} \otimes \Lambda_d \lra \Lambda_i \]  is   
surjective for all $i \le d$. This is so iff the dual map   
\[ \beta_i: (R/I)_i \lra S_{d-i} \otimes (R/I)_d \]   
is injective for all $i$. This map can be written as  
\begin{equation}
\beta_i: u \lra \sum\limits_M \, \, \uy^M \otimes u\, \ux^M,   
\label{expression.beta}   
\end{equation}   
the sum quantified over all monomials $\ux^M$ of degree $d-i$.   

As a consequence, if $R/I$ is level then the graded piece   
$I_d$ determines $I$, by the following recipe:   
$I_i = \{ u \in R_i: u.R_{d-i} \subseteq I_d\}$ for $i \le d$, and   
$I_i = R_i$ for $i > d$. In the terminology of \cite{Ger1} (a related  
terminology was originally introduced by A.~Iarrobino), $I$ is the   
\emph{ancestor ideal} of the vector space $I_d$.   
\end{Remark}   
\begin{Remark} \rm We can detect whether $R/I$ is level from the   
last syzygy module in its minimal resolution. Indeed,  let  
\[ 0 \ra P_n \ra \dots \ra P_0 (=R) \ra R/I \ra 0   
 \qquad (\dagger) \]   
be the graded minimal free resolution of $R/I$, and  
\[ 0 \ra R(-n) \ra R(-n+1)^n \ra \dots \wedge^i (R(-1)^n) \ra \dots   
\ra R \ra k \ra 0 \; \; (\dagger \dagger) \]   
the Koszul resolution of $k$. We will calculate   
the graded $R$-module $N = \text{Tor}_n^R(R/I,k)$ in two ways.   
If we tensor $(\dagger)$ with $k$, then all differentials are zero,   
hence $N = P_n \otimes k$. When we tensor $(\dagger\dagger)$ with   
$R/I$, the kernel in the leftmost place is $N = \text{socle}(A)(-n)$.   
Hence $A$ is level of socle degree $d$ and type $t$ iff $\text{socle}(A)   
\simeq k(-d)^t$, iff $P_n \simeq R(-d-n)^t$.   
\label{remark.socle} \end{Remark}   

Henceforth $A=R/I$ always denotes a level algebra of type $t$ and socle degree   
$d$, loosely said to be of type $(t,d)$.   
Let $\B \subseteq S_d \otimes \O_G$ denote the tautological bundle   
on $G(t,S_d)$, thus its fibre over a point $\Lambda \in G$ is the subspace   
$\Lambda$. The internal products give vector bundle maps   
\begin{equation}
\phi_i: R_{d-i} \otimes \B \lra S_i \otimes \O_G, \quad 1 \le i \le d.   
\label{global.productmap} \end{equation}   
Dually, there are maps   
\[ \phi_i^*: R_i \otimes \O_G \lra S_{d-i} \otimes \B^*. \]   
Now $\B^*$ is the universal quotient   
bundle of $G(\dim R_d -t, \dim R_d)$ via (\ref{dual.grass}),   
thus its fibre over the point $I_d$ is the subspace $R_d/I_d$.   

\subsection{Definition of Level Subschemes}  \label{levelschemes.defn}
We fix $(t,d)$ and let $G = G(t, S_d)$.  
The Hilbert function of $A$ is given by   
\[ H(A,i) = \dim R_i/I_i = \dim \Lambda_i. \]  
This motivates the following definition.  

For integers $i,r$, let  $\Le(i,r)$ be the closed subscheme of $G$ defined  
by the condition $\{ \rank(\phi_i) \le r \}$. (Locally it is defined  
by the vanishing of $(r+1)$-minors of the matrix representing $\phi_i$.)   
Let $\Le^\circ(i,r)$ be the locally closed subscheme   
$\Le(i,r) \setminus \Le(i,r-1)$. Thus $A$ represents a closed point   
of $\Le(i,r)$ (resp. $\Le^\circ(i,r)$) whenever $H(A,i) \le r$   
(resp. $H(A,i) =r$).   

Let $\uh = (h_0, h_1, h_2, \dots)$ be   
a sequence of nonnegative integers such that $h_0 =1, h_d =t$ and $h_i =0$ for   
$i > d$. (It is a useful convention that $h_i =0$ for $i<0$.)  
Define scheme-theoretic intersections
\[ \Le(\uh) = \bigcap_{i=1}^{d-1} \Le(i,h_i), \quad   
   \Le^\circ(\uh) = \bigcap_{i=1}^{d-1} \Le^\circ(i,h_i).   
\]   
These are respectively closed and locally closed subschemes of $G(t,S_d)$.   
Via the identification in (\ref{dual.grass}), we will occasionally think of them as   
subschemes of $G(\dim R_d -t,R_d)$. The point $A=R/I$ lies in $\Le(\uh)$ (resp.   
$\Le^\circ(\uh)$) iff $\dim_k A_i \le h_i$ (resp.~$\dim_k A_i = h_i$)   
for all $i$.   

Of course either of the schemes may be empty, and it is   
in general an open   
problem to characterise those $\uh$ for which they are not. For $n =2$, such   
a characterisation is given in Theorem \ref{theorem.char.level}.   
If $\Le^\circ(\uh)$ is nonempty, then we will say that $\uh$ is a  
\emph{level Hilbert function}.  

\subsection{The structure of $\Le(\uh)$ for $t=1, n=2$.}  
\label{goren.2}  

The structure of the parameter spaces for Gorenstein quotients of $R = k[x_1,x_2]$ is  
rather well-understood and provides a useful paradigm for our study of artin level  
quotients of $R$ having type $>1$.  An outine of this story is given below,  
see \cite{Ger1} for details.  

It is easy to show that $A = R/I$ is a graded  
Gorenstein artin algebra iff $I$ is a complete intersection. Thus $I= (u_1,u_2)$,  
where $u_1,u_2$ are homogeneous and  
$\deg u_1 = a \leq b = \deg u_2$.  In this case $d = a+b-2$, and  
the Hilbert function of $A$ is
\[ H(i) = \begin{cases}
i+1 &      \text{for $0 \le i \le a-1$,} \\  
a   & \text{for $a \le i \le b-1$,} \\  
a+b-(i+1) & \text{for $b \le i \le d$;}
\end{cases} \]  
in particular it is centrally symmetric. We will denote this function by $\uh_a$.  
It follows that for socle degree $d$, there are precisely
\[  
\ell = \begin{cases}
(d+2)/2, & \text{if $d$ is even} \\  
(d+1)/ 2, & \text{if $d$ is odd}  
\end{cases} \]  
possible Hilbert functions for Gorenstein artin quotients of $R$.  
The collection $\{ \uh_a \}$ is totally ordered, i.e.,  
$h_a(j) \leq h_{a+1}(j)$ for all $0 \leq j \leq d$ and $1 \leq a \le \ell$.  
For brevity, let $\Le_a$ denote the scheme  
$\Le(\uh_a) \subseteq \P S_d$.  

In fact, $\Le_a^\circ$ is the locus of \emph{power sums} of length $a$, i.e.,  
\[  
\Le_a^\circ = \{ F \in \P S_d: F = L_1^d + \ldots + L_a^d \quad   
\text{for some $L_i$ in $S_1$} \},
\]  
and $\Le_a$ its Zariski closure. Thus $\Le_1$ can be identified with  
the rational normal curve in $\P S_d$, and $\Le_a$ is the union of  
(possibly degenerate) secant $(a-1)$-planes to $\Le_1$. In particular,  
$\dim \Le_a = 2a-1$.  

Let $z_i = x_1^{d-i} x_2^i$, then $\Sym^\bullet \, R_d = k[z_0, \dots, z_d]$  
is the coordinate ring of $\P S_d$. Consider the Hankel matrix  
\[  
{\mathcal C}_a : =  
\left[ \begin{array}{ccccc}
z_0 & z_1 & \cdots \cdots & z_{d-a} \\  
z_1 & z_2 & \cdots \cdots & z_{d-a+1} \\
\vdots & \vdots & & \vdots  \\
z_a & z_{\ell +1} & \cdots \cdots & z_d  
\end{array} \right]\]  
and let $\wp_a$ denote the ideal of its maximal minors. Then it  
is a theorem of Gruson and Peskine  
(see \cite{GrusonPeskine}) that $\wp_a$ is perfect, prime and equals  
the ideal of $\Le_a$ in $k[z_0, \dots, z_d]$.  
Now the Eagon-Northcott Theorem implies that $\Le_a \subseteq  \P S_d$ is  
an arithmetically Cohen-Macaulay variety.  

\subsection{Summary of results.} In the next section, we derive an expression   
for the tangent space to a point of $\Le^\circ(\uh)$. This is a direct   
generalization of \cite[Theorem 3.9]{IK} to the non-Gorenstein case.   
For sections 3 and 4,   
we assume $n=2$. In section 3, we give a geometric description of a point of   
$\Le^\circ(\uh)$ in terms of secant planes to the rational normal curve, which  
generalises the one given above for $t=1$.  
We relate this description to Waring's problem for systems of algebraic forms   
and solve the   
problem for $n=2$. In the last section we prove a projective normality theorem   
for a class of schemes $\Le(i,r)$ using spectral sequence techniques.   
The results in the following three sections are largely independent   
of each other, and as such may be read separately.   

We thank the referee for several helpful suggestions, and 
specifically for contributing 
Corollary \ref{corollary.count.h}. We owe the result of 
Theorem \ref{theorem.char.level} to G.~Valla. We also 
acknowledge the help of John Stembridge's `SF' Maple package  
for some calculations in \S \ref{catalecticant}.  

\section{Tangent spaces to level subschemes}   
Let $R = k[x_1,\dots, x_n]$ and let $A = R/I$ be an artin level  
quotient of type $(t,d)$. Given  a degree zero  
morphism $\psi: I \lra R/I$ of graded $R$-modules, we have induced maps   
of $k$-vector spaces $\psi_i: I_i \lra (R/I)_i$. We claim that $\psi_d$   
entirely determines $\psi$. Indeed, let $u \in I_i$ and $\ux^M$ a   
monomial of degree $d-i$. Then   
$\psi_d(u \, \ux^M) = \psi_i(u)\,\ux^M$. But then   
$\beta_i(\psi_i(u)) = \sum\limits_M \, \uy^M \otimes \psi_d(u \, \ux^M)$.   
Since $\beta_i$ is injective, this determines $\psi_i(u)$ uniquely.   
Thus we have an inclusion   
\begin{equation} \Hom_R(I,R/I)_0 \hookrightarrow \Hom_k(I_d,R_d/I_d), \quad  
\psi \lra \psi_d.  
\label{incl1} \end{equation}   
We also have a parallel inclusion   
\begin{equation}   
\Hom_R(I^{-1},S/I^{-1})_0 \hookrightarrow \Hom_k(\Lambda,S_d/\Lambda).   
\label{incl2} \end{equation}

Recall that if $U$ is a vector space and $W$ an $m$-dimensional subspace,   
then the tangent space to $G(m,U)$ at $W$ (denoted $T_{G,W}$)   
is canonically isomorphic to $\Hom(W,U/W)$. Thus   
\[   
T_{G(t,S_d),\Lambda} = \Hom_k(\Lambda, S_d/\Lambda) \quad   
T_{G(\dim R_d-t,R_d),I_d} = \Hom_k(I_d, R_d/I_d).   
\]   

\begin{Theorem} Let $A = R/I$ be as above with Hilbert function $\uh$ and  
inclusions (\ref{incl1}), (\ref{incl2}).  
\begin{enumerate}   
\item[(A)]
Regarding $\Le^\circ(\uh)$ as a subscheme of $G(t,S_d)$,   
we have a canonical isomorphism    
\[ T_{\Le^\circ(\uh),\Lambda} = \Hom_R(I^{-1},S/I^{-1})_0 =   
\Hom_R(I^{-1},(I^2)^{-1}/I^{-1})_0. \]
\item[(B)]
Regarding $\Le^\circ(\uh)$ as a subscheme of $G(\dim R_d -t,R_d)$,   
we have a canonical isomorphism    
\[ T_{\Le^\circ(\uh),I_d} = \Hom_R(I,R/I)_0 = \Hom_R(I/I^2,R/I)_0. \]
\end{enumerate}   
\label{theorem.tangentspace} \end{Theorem}   

\demo   
We begin by recalling the relevant result about the tangent space to a   
generic determinantal variety (see \cite[Ch.~2]{ACGH}).   

Let $M = M(p,q)$ denote the space of all $p \times q$ matrices over   
${\mathbf C}$, or   
equivalently the space of vector space maps ${\mathbf C}^{\,p} \ra   
{\mathbf C}^{\,q}$.   
Since $M$ is an affine space, for any $X \in M$, the   
tangent space $T_{M,X}$ can be canonically identified with $M$. Fix   
an integer $r \le   
\min\{p,q\}$ and let $M_r$ be the subvariety of matrices with rank $\le r$. If   
$X \in M_r \setminus M_{r-1}$, then $X$ is a smooth point of $M_r$ and   
\[ T_{M_r,X} = \{ Y \in M: Y (\ker X) \subseteq \image X \}.   
\]   

Now if $\Lambda \in \Le^\circ(i,r)$, then $\phi_i$ is represented   
in a neighbourhood $U \subseteq G(t,S_d)$ of $\Lambda$ by a   
matrix of size $t(d-i+1) \times (i+1)$, whose entries are regular functions on   
$U$. Writing $M = M(t(d-i+1),i+1)$, these functions  
define a morphism $f:U \lra M$.  
Thus the following is a fibre square  
\[ \diagram   
 U \cap \Le^\circ(i,r) \rto \dto & U \dto^f \\   
 M_r \rto_i & M
\enddiagram \]  
Hence   
\[ T_{\Le^\circ(i,r), \Lambda} = \{ \tau \in T_{U,\Lambda} =   
\Hom(\Lambda,S_d/\Lambda): df(\tau) \in T_{M_r,f(\Lambda)} \}. \]   
(Here $df$ denotes the induced map on tangent spaces.)

This translates into the statement that   
$T_{\Le^\circ(i,r),\Lambda}$ consists of all $\tau \in   
\Hom_k(\Lambda, S_d/\Lambda)$ such that   
the broken arrow in the following diagram is zero.

 \[ \diagram
{\ker \, \alpha_i} \rto \drrdashed|>{\tip} &   
R_{d-i} \otimes \Lambda \rto^{{\text{id} \, \otimes \, \tau}}   
& \quad R_{d-i} \otimes S_d/\Lambda \dto^{\mu} \\   
& & S_i/\Lambda_i \enddiagram \]

The map $\mu$ comes from the internal product in an obvious way.   
This implies that $\tau \in T_{\Le^\circ(i,r),\Lambda}$ iff   
the composite $\mu \circ (\text{id} \otimes \tau)$ factors through   
$\image \alpha_i = \Lambda_i$. Let $\tau_i: \Lambda_i \lra S_i/\Lambda_i$   
denote the induced map. Now   
\[ T_{\Le^\circ(\uh),\Lambda} = \bigcap\limits_i \,   
T_{\Le^\circ(i,h_i),\Lambda},   
\]   
hence $\tau \in T_{\Le^\circ(\uh),\Lambda}$ iff it defines a sequence   
$(\tau_i)$ as above which glues to give an $R$-module map   
$I^{-1} \lra S/I^{-1}$. This proves (A).   

For (B), a parallel argument leads to the following: an   
element $\omega \in \Hom_k(I_d, R_d/I_d)$ belongs to   
$T_{\Le^\circ(i,r),I_d}$ iff in the diagram below the broken arrow can be   
filled in.   

\[ \diagram   
I_d \otimes S_{d-i} \rto^{\omega \, \otimes \, \text{id}} & R_d/I_d \otimes S_{d-i} \\   
I_i \uto \rdashed|>{\tip} & R_i/I_i \uto   
\enddiagram \]  

Here both vertical maps are given by formula (\ref{expression.beta}),   
in particular they are injective. Hence the broken arrow is unique if   
it exists, which we then denote by $\omega_i$.   
Thus $\omega \in T_{\Le^\circ(\uh),I_d}$ iff it defines a sequence   
$(\omega_i)$ as above, which glues to give an $R$-module map   
$I \lra R/I$. This proves the theorem. \qed   

\begin{Remark} \rm   
The scheme $\text{GradAlg}(\uh)$ (defined by J.~Kleppe \cite{Kleppe1})   
parametrises graded quotients of $R$ (level or not) with Hilbert function   
$\uh$. Its tangent space at the point $R/I$ is also canonically isomorphic to   
$\Hom(I,R/I)_0$.  
See Remarks 3.10 and 4.3 in \cite{IK} for a more detailed comparison of these  
two spaces (in the Gorenstein case).  
\end{Remark}

\section{Level algebras in codimension two}   
In this section (and the next) we consider quotients of $R=k[x_1,x_2]$.   

\subsection{Preliminaries}
Let $A=R/I$ be an artin level algebra with Hilbert function $H$,   
type $t$ and socle degree $d$.   
By Remark \ref{remark.socle}, we have a resolution    
\begin{equation}
 0 \lra R^{\,t}(-d-2) \lra \bigoplus_{\ell=1}^{d+1} R^{\, e_\ell}(-\ell)  
\lra  R \lra R/I \lra 0.   
\label{ideal.resolution} \end{equation}   
Here $e_\ell$ is the number of minimal generators of $I$ in degree $\ell$, and   
$\sum e_\ell = t+1$. Hence   
\begin{equation}   
H(A,i)  =(i+1)- \sum\limits_{\ell =1}^{i} e_\ell \, (i-\ell+1)   
\quad \text{for all $i \le d+1$}.   
\end{equation}   
With a little manipulation, this implies   
\begin{equation}   
e_{i+1}  = 2H(A,i) - H(A,i-1) - H(A,i+1) \quad \text{for $0 \le i \le d$}.
\label{no.of.gens} \end{equation}   
Hence the sequence $(e_i)$ can be recovered from the Hilbert   
function. Applying the functor $\Hom_R(-,R/I)$ to the resolution of $I$, we   
have an exact sequence   
\[ 0 \lra \Hom(I,R/I) \lra \bigoplus_{\ell=1}^{d+1} (R/I)^{e_\ell}(\ell) \lra   
(R/I)^t(d+2),  \]   
hence   
\begin{equation}
\dim_k \Hom_R(I,R/I)_0 = \sum_{\ell=1}^d e_\ell \,H(A,\ell).   
\label{dimformula} \end{equation}

The next result characterises the level Hilbert functions of type   
$(t,d)$ in codimension two. It is due to G.~Valla, who had kindly communicated 
its proof to the second author a few years ago. A more general version 
(which covers codimension two non-level algebras) is stated by A.~Iarrobino 
in \cite[Theorem 4.6A]{Iarro4}. 
\begin{Theorem}[Iarrobino, Valla] \label{theorem.char.level}   
Let $\uh = (h_0,h_1, \dots)$ be a sequence of nonnegative   
integers satisfying $h_0=1, h_d=t$ and $h_i=0$ for $i > d$.  
Then $\Le^\circ(\uh)$ is nonempty if and only if
\[ 2 h_i \ge h_{i-1} + h_{i+1},  \quad \text{for all \, $0 \le i \le d$.} \]   
(By convention, $h_i = 0$ for $i < 0$.)  
\end{Theorem}   

\demo The `only if' part follows from (\ref{no.of.gens}).   
Assume that $\uh$ satisfies the   
hypotheses. Then we inductively deduce $h_i \le i+1$.  
Define $e_i = 2h_{i-1} - (h_{i-2} + h_i)$ for $1 \le i \le d+1$ and  
$e_i = 0$ elsewhere. Then  
$\sum\limits_{i=1}^{d+1} e_i = h_d + h_0 = t+1$.  
Define a sequence of integers   
\[ {\underline q}: \, q_1 \le q_2 \le \dots \le q_{t+1}, \]
such that for $1 \le i \le d+1$, the integer $i$ occurs $e_i$ times.   
An easy calculation shows that $\sum q_i = \sum i.e_i = t(d+2)$.   

Let $M$ be the $t \times (t+1)$ matrix whose only nonzero entries are   
$M_{i,i}= x_1^{d+2-q_i}$ and $M_{i,i+1}= x_2^{d+2-q_{i+1}}$ for   
$1 \le i \le t$, and let $I$ be the ideal of its maximal minors.   
Since $I$ is $(x_1,x_2)$-primary, it has depth $2$. The $t+1$ maximal  
minors of $M$ are nonzero, and they have degrees $q_1,\dots,q_{t+1}$. By the   
Hilbert-Burch theorem, $R/I$ has a resolution with Betti numbers as in   
(\ref{ideal.resolution}). Then by Remark \ref{remark.socle}, the point   
$A=R/I$ lies in $\Le^\circ(\uh)$. \qed   

\begin{Example} \rm  
Let $(t,d) = (3,7)$ and  
$\uh = (1,\, 2, \, 3, \, 4, \, 5, \, 5, \,  4, \, 3, \, 0)$.  
Then $e_5 = e_6 = 1, e_8=2$, and ${\underline q} = (5,\,6,\,8,\,8)$. Hence  
\[ M = \left[ \begin{array}{cccc}  
x_1^4 & x_2^3 & 0 & 0 \\  
    0 & x_1^3 & x_2 & 0 \\  
    0 &     0 & x_1 & x_2  
\end{array} \right] \]  
and $I = (x_2^5, x_1^4x_2^2, x_1^7x_2, x_1^8)$.  
\label{example.levelfn} \end{Example}  
We owe the following observation to the referee.  
\begin{Corollary} \label{corollary.count.h}
The level Hilbert functions $\uh$ of type $(t,d)$ are in bijection with  
partitions of $d-t+1$ with no part exceeding $t+1$.  
\end{Corollary}  
\demo Given $\uh$, define  
$\mu_i = h_i - h_{i+1} + 1$ for $0 \le i \le d-1$. Then  
${\underline \mu} = (\mu_{d-1}, \dots, \mu_1, \mu_0)$ is a partition as above. Conversely,  
given such a partition we append zeros to make its length equal to $d$,  
and then determine $h_i$ recursively.  
\qed  

\begin{Remark} \rm If $\uh$ is a level Hilbert function, then $\Le^\circ(\uh)$   
is an irreducible and smooth variety. Indeed, the scheme   
$\text{GradAlg}(\uh)$ is irreducible and smooth by a   
result of Iarrobino (\cite[Theorem 2.9]{Ia2}), and   
$\Le^\circ(\uh) \subseteq \text{GradAlg}(\uh)$ is a dense open subset.   
Hence, from (\ref{dimformula}),  
\begin{equation}
 \dim \Le^\circ(\uh) = \sum\limits_{i=1}^d e_ih_i =  
\sum\limits_{i=1}^d \, h_i(2h_{i-1} - h_i -h_{i-2}).
\end{equation}  
E.g., we have $\dim \Le^\circ(\uh) = 9$ in Example \ref{example.levelfn}.  
\end{Remark}   

\subsection{Geometric description of points in $\Le^\circ(\uh)$}  
\label{geom.description}
We start with an example to illustrate the description we have in mind.   
We need the following classical lemma (see \cite[p.~23 ff]{IK},  
also \cite{KungRota}).   
\begin{Lemma}[Jordan] \rm Let $u \in R_m$ be a form factoring as   
\[ \prod\limits_i (a_i x_1 + b_i x_2)^{\mu_i},   
\quad \text{so that $\sum \mu_i = m$.} \]   
If $n \ge m$, then $(u)^{-1}_n$ (the subspace of forms in $S_n$ which   
are apolar to $u$) equals   
\[ \sum_i S_{\mu_i-1}(b_i y_1 - a_i y_2)^{n-\mu_i +1} =   
\{ \sum f_i\, (b_i y_1 -a_i y_2)^{n-\mu_i+1}: f_i \in S_{\mu_i-1} \}.
\]
In particular, this is an $m$-dimensional vector space.  
\end{Lemma}  \qed   

\begin{Example} \rm   
Let $(t,d) = (2,6)$ and consider the level Hilbert function   
\[ \uh=(1, \, 2, \, 3, \, 4, \, 4, \, 3, \, 2, \, 0). \]   
Then $\Lambda \in \Le^\circ(\uh)$ defines a line $\P \Lambda$ in  
$\P^6 (= \P S_6)$. We identify the subset   
$C_6 = \{ [L^6] \in \P S_6 : L \in S_1 \}$ as the rational normal sextic in   
$\P S_6$.   

By formula (\ref{no.of.gens}), $I = \ann(\Lambda)$ has one minimal   
generator each in degrees $4,5,7$. Let $u_4 \in R_4$ be the first generator,    
factoring as $u_4 = \prod\limits_{i=1}^4 (a_i x_1 + b_i x_2)$.   
 For simplicity, assume that    
$[a_1,b_1], \dots, [a_4,b_4]$ are distinct points in $\P^1$. Then   
by Jordan's lemma, the subspace $(u_4)^{-1}_6 \subseteq S_6$   
is the span of $(b_iy_1 - a_i y_2)^6, 1 \le i \le 4$. Let   
$\Pi_4$ denote the projectivisation $\P (u_4)^{-1}_6$,   
which is the secant $3$-plane to $C_6$   
spanned by the four points $[b_i,-a_i]$. Consider generators $u_5,u_7$ and   
define $\Pi_5, \Pi_7$ analogously. (Of course, $\Pi_7 = \P S_6$.)   
Now   
\[ (I^{-1})_6 = ((u_4,u_5,u_7)^{-1})_6 \implies   
\P \Lambda = \Pi_4 \cap \Pi_5 \cap \Pi_7. \]   
Thus (the line corresponding to)  
every element $\Lambda \in \Le^\circ(\uh)$ is representable  
as an intersection of secant planes to the rational normal curve,    
in a way which depends only on the combinatorics of $\uh$.   
If (say) $u_4$ has multiple roots, then $\Pi_4$ is tangent to the curve at   
one or more points, so must be counted as a degenerate secant plane.   
\end{Example}   

\begin{Definition} \rm   
Let $C_d = \{ [L^d]: L \in S_1 \}$ be the rational normal curve in $\P S_d$.   
A linear subspace $\Pi \subseteq \P S_d$ of (projective) dimension $s$    
will be called a \emph{secant $s$-plane to $C_d$} if
the scheme-theoretic intersection $C_d \cap \Pi$ has length $ \ge s+1$.   
(Then the length must equal $s+1$, essentially by Jordan's lemma.)
\end{Definition}
Now for an arbitrary level algebra in codimension two, we have the following   
description.   
\begin{Proposition}
Let $\uh$ be a level Hilbert function and $\Lambda \in \Le^\circ(\uh)$.   
Define a sequence ${\underline q}$ as in the proof of   
Theorem \ref{theorem.char.level}. Then there exist secant $(q_i-1)$-planes   
$\Pi_{q_i}$ such that   
\begin{equation}
\P \Lambda = \Pi_{q_1} \cap \dots \cap \Pi_{q_{t+1}}.   
\label{lambda.expression} \end{equation}
\label{prop.sec.1} \end{Proposition}
\demo The essential point already occurs in the example above.   
Let $q$ be one of the $q_i$ and $u \in I_q$ a generator.   
Write   
\[ u = \prod\limits_i (a_i x_1 + b_i x_2)^{\mu_i}. \]
Let $\Phi_i$ be the osculating $(\mu_i-1)$-plane to $C_d$ at   
the point $(b_i y_1 - a_i y_2)^d$. (This is the point itself if   
$\mu_i =1$.) Algebraically, $\Phi_i$ is the projectivisation   
$\P (S_{\mu_i-1} (b_i y_1 - a_i y_2)^{d- \mu_i+1})$.   

Let $\Pi_q$ be the linear   
span of all the $\Phi_i$. Then $\dim \Pi_q=q-1$ and by Jordan's lemma,   
$\Pi_q$ is the locus of forms $F \in S_d$ apolar to $u$. Construct such a   
plane $\Pi_{q_i}$ for each $q_i$. A form $F$ lies in   
$\Lambda$ iff it is apolar to each generator of $I$, iff it belongs to   
$\bigcap \Pi_{q_i}$. The proposition is proved. \qed   

\begin{Remark}  \rm   
The argument heavily depends on the fact that any zero-dimensional  
subscheme of $C_d$ is in linearly general position. This property characterises  
rational normal curves (see \cite[p.~270]{GH}).
\end{Remark}   
   
The preceding proposition admits a converse. Consider the following example.   
\begin{Example} \rm   
Let $\P \Lambda \subseteq \P S_{11} ( = \P^{11})$ be a line appearing  
as an intersection   
\[ \P \Lambda = \Pi_8 \cap \Pi_8' \cap \Pi_{10},   
\]
where each $\Pi_q$ is a secant $(q-1)$-plane to $C_{11}$. Note that the   
planes intersect properly, i.e., in the expected codimension. We claim that    
$\Lambda$ belongs to $\Le^\circ(\uh)$, for   
\[ \begin{array}{ccccccccccccc}
\uh = 1 & 2 & 3 & 4 & 5 & 6 & 7 & 8 & 7 & 6 & 4 & 2 & 0.
\end{array} \]   
Let $W_8 \subseteq S_{11}$ be the subspace such that $\Pi_8 = \P W_8$ etc.   
Assume $\Pi_8$ intersects $C_{11}$ at points $(b_i y_1 - a_i y_2)^{11}$   
with respective multiplicities $\mu_i$ (so that $\sum \mu_i = 8$). Let   
$u_8 \in R_8$ be the element $\prod (a_i x_1 + b_i x_2)^{\mu_i}$, then   
$\ann(W_8)$ is the principal ideal $(u_8)$. Define elements $u'_8,u_{10}$   
similarly. Now   
\[ I = \ann(\Lambda) = \ann(W_8) + \ann(W_8') + \ann(W_{10}) =   
(u_8,u_8',u_{10}).   
\]   
Since by construction, $R/I$ is a level algebra of type two, $I$  
has three minimal   
generators. Hence their degrees must be $(8,8,10)$ and then the Hilbert function   
of $R/I$ is determined by formula (\ref{no.of.gens}). This completes   
the argument.   
\end{Example}   

\begin{Proposition}   
Let $\P \Lambda \subseteq \P S_d$ be a $(t-1)$-dimensional subspace which   
is expressible as an intersection   
\[ \P \Lambda = \Pi_{q_1} \cap \dots \cap \Pi_{q_{t+1}}, \]
where \begin{itemize}   
\item   
$\Pi_{q_i}$ is a secant $(q_i-1)$-plane to $C_d$,   
\item   
the intersection is proper, i.e.,   
\[ \sum\limits_i \text{codim}\,(\Pi_{q_i}, \P^d) =  
\text{codim}\,(\P \Lambda, \P^d) = d-t+1. \]  
\end{itemize}   
Then $R/\ann(\Lambda) \in \Le^\circ(\uh)$
where $h_0=1$ and $h_i = 2h_{i-1} - h_{i-2} - e_i$ for   
$1 \le i \le d$. The $e_i$ are defined by arranging the $q_i$ as  
in Theorem \ref{theorem.char.level}.  
\label{prop.sec.2} \end{Proposition}   

\noindent \demo Left to the reader.  \qed   

Propositions \ref{prop.sec.1}, \ref{prop.sec.2} are natural generalizations  
of the description of $\Le(\uh)$ in the Gorenstein case.  

\begin{Example} \rm The minimal Hilbert function of type $(t,d)$ is   
$\uh = (1,\, 2, \dots, t-1,\, t, \dots, t, 0)$ and then   
$\Le^\circ(\uh) = \Le(\uh)$ is   
the variety of secant $(t-1)$-planes to $C_d$. Abstractly,   
$\Le(\uh) \simeq  \Sym^t \, \P^1 \simeq \P^t$.   
See \cite{BiGer} for a description of the minimal level Hilbert function  
when $n > 2$.  

The maximal function is $h_i = \min \{ i+1, (d-i+1)t \}$ for $1 \le i \le d$.  
If we let  
\[ s_0 = \lceil \frac{t(d+1)}{t+1} \rceil,  \]  
then $h_{s_0} < s_0 +1$. Hence every $\Lambda \in G(t,S_d)$ has an apolar  
form of degree $\le s_0$.  
\label{example.minmax} \end{Example}   

Using this formulation, the so-called Waring's problem can be solved   
completely for binary forms.   
\subsection{Waring's problem for several binary forms}   
\label{section.waring}   
We will start with an informal account.   
Given binary forms $F_1,\dots,F_t$ of degree $d$, we would like to find   
linear forms $L_1,\dots,L_s$ such that   
\begin{equation}   
 F_i = c_{i1} L_1^d + \dots + c_{is} L_s^d,   
\label{add.expression} \end{equation}   
for some $c_{ij} \in k$. This is always possible for $s=d+1$,  
indeed if we choose $L_1, \dots, L_{d+1}$ generally, then  
the $L_i^d$ span $S_d$.   

The \lq simultaneous\rq \, Waring's problem, in one of its versions,  
is to find the smallest $s$ which suffices   
for a \emph{general} choice of $F_i$. (See Bronowski \cite{Bronowski1}  
for a discussion of $n$-ary forms.) Here we consider a more general   
version, i.e., we fix $s$ and consider the   
locus $\Sigma$ of forms $\{F_1,\dots,F_t\}$ which admit such a representation.   
In practice, we allow not only representations as above, but also   
generalised additive decompositions (GADs) in the sense  
of \cite[Definition 1.30]{IK}.   

\begin{Definition} \rm   
A finite collection $\uL = \{L_i\} \subseteq S_1 \setminus \{0\}$   
of linear forms will be called \emph{admissible} if any  
two $L_i,L_j$ are nonproportional.   
\end{Definition}   

\begin{Definition} \rm   
Let $\uL$ be an admissible collection as above and $F \in S_d$. A  
\emph{GAD for  $F$ with respect to  $\uL$} is a  
collection $\{p_i\}$ of forms such that   
$F = \sum\limits_i p_i \, L_i^{d - \alpha_i}$, and  
$\alpha_i = \deg p_i < d$. (We set $\deg \, 0 = -1$ by convention.)
The integer $\ell_{\ualpha} = \sum\limits_i (\alpha_i + 1)$ is called  
the \emph{length} of the GAD.  
\end{Definition}   
If all $\alpha_i =0$, then this reduces to the expression (\ref{add.expression}).  
Given sequences $\ualpha$ and $\uL$ such that $-1 \le \alpha_i < d$,  
define a subspace  
\[ W(\ualpha, \uL) = \sum S_{\alpha_i} L^{d-\alpha_i} \subseteq S_d.  \]   
(By convention, $S_{-1} = 0$.)
Write $L_i = b_i y_1 - a_i y_2$, then by hypothesis,  
$[a_i,b_i]$ represent distinct points of $\P^1$.   
\begin{Lemma}[\cite{IK}, \S 1.3]   
With notation as above,   
\begin{itemize}   
\item   
A form lies in $W(\ualpha, \uL)$ iff it is apolar to   
$\prod\limits_i (a_i x_1 + b_i x_2)^{\alpha_i+1}$.   
In particular, $\dim W(\ualpha,\uL) = \sum\limits_i (\alpha_i +1)$.   
\item   
Let $\Phi_{i,\alpha_i}$ be the osculating $\alpha_i$-plane to $C_d$   
at the point $L_i^d \in C_d$. (By convention empty if $\alpha_i = -1$.)   
Then the linear span of all $\{\Phi_{i,\alpha_i}\}_i $ is  
the projectivisation $\P W(\ualpha, \uL)$.  
\end{itemize} \end{Lemma}   
This is merely a rephrasing of Jordan's lemma.   
Thus the subspaces $W(\ualpha, \uL)$ exactly correspond to   
secant $(\ell_\ualpha-1)$-planes to $C_d$.   
Now let $(d,t)$ be as above and fix an integer $s$ such that   
$d+1 > s \ge t$. Define   
\begin{equation}    
\begin{aligned}   
\Sigma_s =   & \{ \Lambda \in G(t,S_d): \,   
\P \Lambda \; \text{lies on some secant $(s-1)$-plane of $C_d$}\}. \\   
\end{aligned} \end{equation}   
Algebraically, $\Lambda \in \Sigma_s$ iff there exists an admissible collection   
$\{L_i\}$ such that each $F \in \Lambda$ has a GAD of length $\le s$ with   
respect to $\{L_i\}$. There is no loss of generality in assuming that   
$\{L_i\}$ has cardinality $s$.   

Now Waring's problem can be interpreted as one of calculating $\dim \Sigma_s$.   
Evidently, it is bounded by $\dim G(t,S_d) = t(d-t+1)$. Let   
$U_s \subseteq \Sym^s (\P S_1)$ be the open subset of admissible   
collections $\uL = \{ L_1, \dots, L_s \}$, and consider the incidence   
correspondence   
\[ \begin{aligned}   
{\widetilde \Sigma_s} =  \{ (\Lambda, \uL) & \in G(t,S_d) \times U_s:   
\; \text{each element of $\Lambda$} \\   
& \text{admits a GAD of length $\le s$ w.r.t. $\uL$} \}.   
\end{aligned} \]   
Then $\Sigma_s$ is the image of the projection   
$\pi_1: {\widetilde \Sigma_s} \lra G(t,S_d)$.   

\begin{Lemma} Each fibre of the projection   
$\pi_2: {\widetilde \Sigma_s} \lra U_s $ is of dimension $t(s-t)$,   
hence $\dim {\widetilde \Sigma_s} =s+t(s-t)$.   
\end{Lemma}   
\demo Each $\Lambda \in \pi_2^{-1}(\uL)$ is a subspace of the (finite) union   
\[ \bigcup\limits_{\ell_{\ualpha} \leqslant s} W(\ualpha,\uL), \]   
hence it is a subspace of one of the $W$. Now use the fact that   
$\dim G(t,W)$ equals $t (\dim W -t)$. \qed   

Thus we have a na{\" \i}ve estimate   
\begin{equation} \dim \, \Sigma_s \le   
\min\, \{\, \underbrace{s+t(s-t)}_{N_1},   
\,\underbrace{t(d-t+1)}_{N_2} \, \}.  
\label{dimsigma.s} \end{equation}

\begin{Theorem}
In (\ref{dimsigma.s}), we have an equality.  
\label{theorem.dimsigma} \end{Theorem}
\demo   
Assume $N_1 < N_2$. We will exhibit a level Hilbert function $\uh$ of type   
$(t,d)$ such that $h_s = s$ and $\dim \Le^\circ(\uh) = N_1$.  
Then by construction, for each $\Lambda \in \Le^\circ(\uh)$ there is  
a nonzero form in $R_s$ which is apolar $\Lambda$.  
This form defines a secant $(s-1)$-plane containing $\P \Lambda$.   
Hence $\Le^\circ(\uh) \subseteq \Sigma_s$, which forces  
$\dim \Sigma_s = N_1$.   

Let $m$ be the unique integer such that $(m+1)t > s \ge mt$.   
Then $N_1 < N_2$ forces $s \le d-m$. Define a sequence $\uh$ by   
\[ h_i = \begin{cases}
i+1 & \text{for $0 \le i \le s-1$,} \\
s   & \text{for $s \le i \le d-m$,} \\
(d-i+1)t & \text{for $d-m+1 \le i \le d$,} \\
0 & \text{for $i > d$.}
\end{cases} \]
It is an immediate verification that $\uh$ is level. Now  
\[ e_s=1, \; e_{d-m+1} = s-mt, \; e_{d-m+2} = mt+t-s,  
\] and $e_i=0$ elsewhere, hence $\dim \Le^\circ(\uh) = N_1$.  

E.g., for $(t,d,s) = (3,11,7)$,   
\[ \begin{array}{ccccccccccccc}
\uh =  1 & 2 & 3 & 4 & 5 & 6 & 7 & 7 & 7 & 7 & 6 & 3 & 0.
\end{array} \]   

If $N_1 \ge N_2$, then $s \ge s_0$ (defined as in Example \ref{example.minmax}).  
But since every $\P \Lambda$ lies on an $(s_0-1)$-secant,  
$\Sigma_{s_0} = G(t,S_d)$ and we are done.  
\qed

\smallskip  

The theorem implies that, given \emph{general} binary forms $F_1, \dots, F_t$, a  
reduction to the expression (\ref{add.expression}) is always possible, as long as  
we have sufficient number of constants implicitly available on the right hand side of  
(\ref{add.expression}). This is no longer so for $n \ge 3$, and the corresponding  
reduction problem is open. See \cite{ego.car} for an approach, where  
Theorem \ref{theorem.dimsigma} is proved using a different method.  

A length $s$-subscheme of $C_d$ is called a polar $s$-hedron  
of $\Lambda$, if $\Lambda$ lies on the corresponding $(s-1)$-secant plane.  
We have shown that the variety of polar $s$-hedra of   
$\Lambda$ is the projective space $\P(R_s/\ann(\Lambda)_s)$.   
For $n \ge 3$, the geometry of this variety is rather more   
mysterious--see \cite{DolgachevKanev, RanestadSch}.   

\subsection{An analogue of the catalecticant}   
\label{catalecticant}
An interesting special case occurs when $N_1 = N_2 -1$, i.e.~when   
$\Sigma_s$ is a hypersurface in $G(t,S_d)$.   
This is  possible iff   
\[ (t+1) | \, (d+2), \; \text{and}\;  s = d+1 - \frac{d+2}{t+1}. \]   
Then $\Sigma_s$ is set-theoretically equal to $\Le(s,s) =   
\{ \rank (\B \otimes R_{d-s} \lra S_s) \le s \}$. Now $\Le(s,s)$ is  
the zero scheme of a global section of the line bundle   
\[ \wedge^{s+1}(\B^* \otimes S_{d-s}) = \O_G(d-s+1). \]   
Hence the fundamental class $[\Le(s,s)]$ equals  
$(d-s+1)\, c_1(\B^*) \in H^2(G,\integers)$.   

The rank condition above can be written as a determinant.   
For instance, let $t=2,d=7,s=5$, and let   
\[ F_1 = \sum\limits_{i=0}^7 a_i y_0^{7-i} y_1^i, \quad   
   F_2 = \sum\limits_{i=0}^7 b_i y_0^{7-i} y_1^i \]   
be linearly independent forms. Then the pencil $\P \Lambda = \P (\vspan(F_1,F_2))$   
lies on a secant $4$-plane to $C_7$ iff   
\[ \left| \begin{array}{cccccc}   
a_0 & a_1 & a_2 & a_3 & a_4 & a_5 \\   
a_1 & a_2 & a_3 & a_4 & a_5 & a_6 \\   
a_2 & a_3 & a_4 & a_5 & a_6 & a_7 \\   
b_0 & b_1 & b_2 & b_3 & b_4 & b_5 \\   
b_1 & b_2 & b_3 & b_4 & b_5 & b_6 \\   
b_2 & b_3 & b_4 & b_5 & b_6 & b_7 \end{array} \right| =0. \]   
This is the analogue of the catalecticant for systems of binary forms.   

\begin{Example} \rm   
Even if $\Sigma_s$ is not a hypersurface, such determinantal conditions   
can always be written down. E.g., let $t=2,d=5, \dim G(2,S_5)=8$.   
The possible level Hilbert functions are   
\[ \begin{array}{cccccccc}   
\uh_1: & 1 & 2 & 2 & 2 & 2 & 2 & 0 \\   
\uh_2: & 1 & 2 & 3 & 3 & 3 & 2 & 0 \\   
\uh_3: & 1 & 2 & 3 & 4 & 3 & 2 & 0 \\   
\uh_4: & 1 & 2 & 3 & 4 & 4 & 2 & 0   
\end{array} \]   
with $\dim \Le^\circ(\uh) = 2,5,6,8$ respectively.   
Now $\Sigma_2 = \Le(\uh_1)$, which is set-theoretically equal to   
any of the schemes $\Le(i,2), i = 2,3,4$ (as defined in \S \ref{levelschemes.defn}).  
Similarly $\Sigma_3 = \Le(\uh_2)$ which is set-theoretically $\Le(3,3)$.   
It is clear that we can write the condition for $\Lambda \in G$ to lie in   
$\Sigma_2$ (or $\Sigma_3$) as the vanishing of certain minors. Indeed, in   
general this can be done in more than one way.   

We can calculate the classes of these schemes  
in $G(t,S_d)$ by the Porteous formula.   
We explain this briefly, see \cite[ \S 14.4]{Fu1} for the details.   

Let $E \stackrel{\alpha}{\lra} F$ be a morphism of vector bundles on a   
smooth projective variety. Assume that $E, F$ have ranks $e,f$ respectively,   
and that the locus $X_r = \{ \rank \, \alpha \le r \}$ is   
of pure codimension   
$(e-r)(f-r)$ in the ambient variety. Then the fundamental class   
of $X_r$ (which we denote by $[X_r]$)
equals the $(e -r) \times (e-r)$  determinant     
whose $(i,j)$-th entry equals the $(f -r + j -i)$-th Chern class   
of the virtual bundle $F - E$. By the Whitney product formula,   
the total Chern class $c_t(F - E) = c_t(F)/c_t(E)$.   

We will follow the conventions of \cite[\S 14.7]{Fu1} for Schubert calculus.
Thus the $i$-th Chern class of the tautological bundle   
$\B$ is $(-1)^i \{1, \dots, 1\}$, where $1$ occurs $i$ times.   
For $G(2, S_5)$, we have $c_t(\B) = 1 - \{1\} + \{1,1\}$.   
A straightforward calculation (done using   
the Maple package `SF') shows that   
\[ [\Le(2,2)] = [\Le(4,2)] = 10\{3,3\} + 6 \{4,2\}. \]   
The formula does not apply to   
$\Le(3,2)$, since it fails to satisfy the codimension hypothesis.   
By a similar calculation, $\Le(3,3) = 8 \{2,1\}$.   

For any $I_5 \in \Le(\uh_1)$, a map in $\Hom_R(I,R/I)$ is entirely   
determined by the image of the unique generator in $I_2$. Now the   
proof of Theorem \ref{theorem.tangentspace} shows that   
the space $T_{\Le(2,2),I_5}$ must be $2$-dimensional, which implies   
$\Le(2,2) = \Le(\uh_1)$.   
By the Littlewood-Richardson rule,   
\[ [\Le(\uh_1)].\{1,1\} = 10 \{4,4\}, \]   
i.e., a general hyperplane $\Omega \subseteq \P S_5$ contains ten secant   
lines of the rational normal quintic $C_5$. Of course these are   
the pairwise joins of the five points $\Omega \cap C_5$. Similarly,   
\[ [\Le(\uh_1)].\{2,0\} = 6 \{4,4\}, \]   
i.e., there are six secant lines to $C_5$ touching a general $2$-plane   
$\Psi$. This can be seen differently: the projection from $\Psi$   
maps $C_5$ onto a rational nodal quintic in $\P^2$, and the six secants   
give rise to the six nodes of the image.  
\end{Example}   

\section{Free resolutions of level subschemes}   
\label{free.res}
We continue to assume $n=2$. Since the scheme $\Le(i,r)$  
is a degeneracy locus in the sense of \cite[Ch.~2]{ACGH}, we can describe its   
minimal resolution following Lascoux \cite{La}, provided it has   
the `correct' codimension in $G(t,S_d)$.   

This granted, in the presence of an additional   
numerical hypothesis (explained below), we can deduce   
that it is arithmetically Cohen-Macaulay in the Pl{\"u}cker embedding.   
In particular we get another proof of the known fact   
that $\Le(i,r)$ is always ACM in the Gorenstein case. In the sequel,   
we need the Borel-Weil-Bott theorem for calculating the cohomology   
of homogeneous vector bundles on Grassmannians. We refer to   
\cite{ego2} for an explanation of the combinatorics involved, but   
also see \cite[p.~687]{Porras}.   

Recall the definition of $\Le(i,r)$ given in \S \ref{levelschemes.defn}. To avoid  
trivialities, we assume $ r < i+1$ throughout the section.  
From \cite[Ch. 2]{ACGH}, we have the following estimate:   
if $c$ is the codimension of any component of $\Le(i,r)$ in $G(t,S_d)$, then   
\begin{equation} \begin{aligned}   
c \le & \; (\rank (\B \otimes R_{d-i}) - r)(\rank \, S_i - r) \\   
= & \; (t (d-i+1) -r)(i+1-r). \end{aligned}   
\label{dim.estimate} \end{equation}
Consider the following conditions:   
\begin{enumerate}   
\item[(C1)] the scheme $\Le(i,r)$ is equidimensional and equality holds in   
(\ref{dim.estimate}) for each component;  
\item[(C2)] $r - (d-i)(i-r) \ge t$.   
\end{enumerate}   

We will impose conditions C1,C2 on the data $(t,d,i,r)$. The next result  
shows the rationale behind C2, as well as its scope of validity.  

\begin{Lemma} \label{lemma.c2}  
The condition C2 holds iff there is a level algebra  
$A$ of type $(t,d)$   
such that $H(A,i-1) = i$ and $H(A,i) \le r$.  
\end{Lemma}   
\demo   
Given the existence of $A$, we have  
\[ \begin{aligned} r - t  & = H(A,i) - H(A,d) =   
\sum\limits_{j=i}^{d-1} H(A,j) - H(A,j+1) \\   
& \ge (d-i)(H(A,i-1) - H(A,i)) \ge (d-i)(i-r),   
\end{aligned} \]   
where the first inequality follows from Theorem \ref{theorem.char.level}.   
Conversely, assume C2 and define  
\[ h_j = \begin{cases}  
j+1 & \text{for $0 \le j \le i-1$,} \\  
\min \{ r - (j-i)(i-r), t(d-j+1) \} & \text{for $i \le j \le d$,} \\
0 & \text{for $j > d$.}  
\end{cases} \]  
Then $\uh$ satisfies the hypotheses of Theorem \ref{theorem.char.level},  
hence it is the Hilbert function of a level algebra $A$. Evidently,  
$h_{i-1} = i$ and $h_i \le r$.  
\qed  

Now assume that $\Le = \Le(i,r)$ satisfies C1 (but not necessarily C2)   
and let $c = \text{codim} \, \Le$.   
Then by the central result of \cite{La}, we have a locally  
free resolution:   
\begin{equation} \begin{aligned}
0 \ra \E^{- c} \ra \dots \ra   
\E^p \ra \E^{p+1} \ra \dots \ra \E^0   
& (= \O_G) \ra \O_{\Le} \ra 0, \\
& \text{for $-c \le p \le 0$;}   
\end{aligned} \label{Lascoux.res} \end{equation}
where   
\begin{equation} \E^p = \bigoplus\limits_{\nu(\lambda')-|\lambda|=p}
\Sc_\lambda (\B \otimes R_{d-i}) \otimes   
H^{\nu(\lambda')}(G', \Sc_{\lambda'} \,Q^*_{G'}).
\end{equation}
This is to be read as follows:   
$G'$ denotes the Grassmannian $G(r,S_i)$, and $Q_{G'}$ its universal quotient   
bundle. The $\lambda$ denote partitions, and   
$\Sc_{\lambda}$ the corresponding Schur functors (where we follow the   
indexing conventions of \cite[Ch.~6]{FH}).   

Given a partition $\lambda$, the Borel-Weil-Bott theorem implies that the   
bundle  $\Sc_{\lambda'} \,Q^*_{G'}$ (resident on $G'$)   
has nonzero cohomology in at most   
one dimension. This number (if it exists) is labelled $\nu(\lambda')$. The direct sum is   
quantified over all $\lambda$ such that $\nu(\lambda')$ is defined. Since   
$\lambda$ has at most $t(d-i + 1)$ rows and $i-r+1$ columns, the sum is finite.

\begin{Example} \rm   
Let $(t,d,i,r) = (2,7,5,4)$ and consider the level Hilbert functions   
\[ \begin{array}{cccccccccc}
\uh_1: & 1 & 2 & 3 & 4 & 5 & 4 & 3 & 2 & 0  \\   
\uh_2: & 1 & 2 & 3 & 4 & 4 & 4 & 4 & 2 & 0   
\end{array} \]   
Then $\Le = \Le(5,4)$ is a union of two components   
$\Le(\uh_1), \Le(\uh_2)$, each of codimension $4$. Hence C1 is satisfied, and   
we have a resolution   
\[   
0 \ra \E^{-4} \ra \dots \ra \E^{-1} \ra \O_{G(2,S_7)} \ra \O_{\Le} \ra 0,   
\]   
where   
\[ \begin{array}{l}   
\E^{-1} = \wedge^5(\B \otimes R_2) \otimes S_5, \\   
\E^{-2} = \wedge^6(\B \otimes R_2) \otimes \Sc_{(21111)}(S_5)   
\oplus \Sc_{21111} ( B \otimes R_2), \\   
\E^{-3} = (\B \otimes R_2) \otimes \wedge^6(\B \otimes R_2)   
\otimes \wedge^5 S_5, \\   
\E^{-4} = [\wedge^6(\B \otimes R_2)]^{\otimes 2}. \\   
\end{array} \]   
The ranks of $\E^0, \dots, \E^{-4}$ are $1,36,70,36,1$, hence $\Le$ is a   
Gorenstein scheme.   
This resolution is equivariant with respect to the action of $SL_2$ on   
the embedding $\Le \subseteq G$.   
\label{example.lascoux} \end{Example}   

The term $\Sc_\lambda (\B \otimes R_{d-i})$ decomposes as a direct sum
\begin{equation} \bigoplus\limits_{\rho,\mu}\; (\Sc_\rho \, \B \otimes \Sc_\mu    
R_{d-i})^{C_{\lambda\rho\mu}},   
\label{kronecker} \end{equation}
quantified over all partitions $\rho,\mu$ of $|\lambda|$.   
The coefficients $C_{\lambda\rho\mu}$ come from the Kronecker   
product of characters of the symmetric group. We explain this briefly,   
see \cite[p.~61]{FH} for details. Also see \cite{Coleman} for   
a tabulation of $C_{\lambda\rho\mu}$ for small values of $|\lambda|$.   

Let $\lambda, \rho, \mu$ be partitions of an integer $a$, and let   
$R_\lambda, R_\rho, R_\mu$ denote the corresponding irreducible   
representations of the symmetric group on $a$ letters (in characteristic   
zero). Then $C_{\lambda \rho \mu}$ is the number of trivial  
representations in the tensor product  
$R_{\lambda} \otimes R_{\rho} \otimes R_{\mu}$. In particular,  
this number is symmetric in the three partitions involved. The  
main combinatorial result that we need is a direct corollary  
of \cite[Theorem 1.6]{Dvir1}.  
\begin{Theorem}[Dvir]
With notation as above, assume $C_{\lambda \rho \mu} \neq 0$. Then   
\[ \begin{aligned}  
\rho_1 & \; (\text{the largest part in $\rho$})  \\
\le & \; (\text{number of parts in $\lambda$}) \cdot (\text{number of parts in $\mu$}).
\end{aligned} \]  
\end{Theorem}

Now we come to the main theorem of this section.   
Recall that a closed subscheme $X \subseteq \P^N$ is said  
to be projectively normal, if the map  
$H^0(\P^N,\O_{\P}(m)) \lra H^0(X,\O_X(m))$ is  
surjective for $m \ge 0$.   

We will regard $\Le(i,r)$ as a closed subscheme of   
$\P (\bigwedge^t S_d)$ via the Pl{\"u}cker embedding of $G(t,S_d)$.   
\begin{Theorem} \label{acm}
Assume that the data $(t,d,i,r)$ satisfy C1.   
\begin{itemize}   
\item[(a).]   
If C2 holds, then $\Le(i,r)$ is projectively normal.   
\item[(b).]
Moreover, if either $t=1$ or C2 is a strict inequality,   
then $\Le(i,r)$ is arithmetically Cohen-Macaulay.   
\end{itemize}   
\end{Theorem}
\demo We will use the following criterion (see \cite[p.~467]{Ei}):
an equidimensional closed subscheme $X \subseteq \P^N$   
(of $\dim > 0$)  
is arithmetically Cohen-Macaulay (ACM) iff   
it is projectively normal and   
$H^j(X, \O_X(m))=0$ for all $m \in {\mathbf Z}$   
and $0 < j < \dim\, X$.  
\smallskip

Since the Grassmannian is projectively normal in the   
Pl{\"u}cker embedding (see e.g.~\cite{ACGH}), part (a)   
will follow if the map   
$H^0(\O_G(m)) \lra  H^0(\O_\Le(m))$ is shown to be surjective   
for $m \ge 0$.   
\medskip

For $m \in {\mathbf Z}$, we have a hypercohomology spectral sequence   
coming from the resolution (\ref{Lascoux.res}).   
\begin{equation} \begin{aligned}   
{} & E_1^{p,q} = H^q(G(t,S_d), \E^p(m)), \;\; d_r^{\,p,q} \lra d_r^{\,p+r,q-r+1}, \\
   & E_\infty^{p,q} \RA H^{p+q}(\O_\Le(m)).   
\end{aligned} \end{equation}
The terms live in the second quadrant, specifically in the   
range   
\[ - c \le p \le 0, \; 0 \le q \le t(d-t+1).   
\]   
Now the theorem will follow from the following lemma:

\begin{Lemma} \begin{enumerate}
\item   
Assume that C2 holds. Then $E_1^{p,q} =0$ for   
$m \ge 0, \, q \neq 0$.   
\item   
Assume that either $t=1$ or C2 is strict. Then   
$E_1^{p,q} =0$ for $m < 0$ and $q \neq \dim G(t,S_d)$.   
\end{enumerate}   
\label{acm.lemma} \end{Lemma}

Let us show that the lemma implies the theorem. Firstly assume $m \ge 0$   
and C2 holds. Then the only nonzero term on the diagonal   
$p+q=0$ is at $p=q=0$. Hence   
$E_\infty^{0,0}=H^0(\O_\Le(m))$ is a quotient of   
$E_1^{0,0} = H^0(\O_G(m))$, which proves (a).   

Now assume $m$ arbitrary, and that either C2 is strict or $t=1$.   
Let $(p,q)$ be such that   
$0 < p+q < \dim \Le$. Then $E_1^{p,q} = 0$, which implies   
$H^j(\O_{\Le}(m)) =0 $ for $j \neq 0, \dim \Le$. This proves (b). \qed   
\medskip   

\noindent {\sc Proof of Lemma \ref{acm.lemma}.}\,
Let $p,q$ be such that $E_1^{p,q} \neq 0$. By hypothesis, $\E^p$ has a summand   
\[ {\mathcal A} = \Sc_\rho \,\B \otimes \Sc_\mu    
R_{d-i} \otimes H^{\nu(\lambda')}(G', \Sc_{\lambda'}\, Q^*_{G'})\]
such that $H^q(G(t,S_d), {\mathcal A}(m)) \neq 0$. Now,
\begin{itemize}
\item{ $\Sc_\mu R_{d-i} \neq 0$ implies that $\mu$ has at most $d-i+1$ rows.}
\item{$\Sc_{\lambda'}\, Q^*_{G'} \neq 0$ implies that $\lambda'$ has at most   
$i-r+1$ rows, i.e., $\lambda$ has at most $i-r+1$ columns. }
\end{itemize}
But then by Dvir's theorem, $C_{\lambda\rho\mu} \neq 0$ implies   
$\rho_1 \le (d-i+1)(i-r+1)$.   
The next step is to use the Borel-Weil-Bott theorem on $\Sc_\rho\,\B   
\otimes \O_G(m)$.   
Let $\gamma$ be the sequence   
$(\underbrace{m,\dots,m}_{d-t+1}; \rho_1,\dots,\rho_t)$. Since   
$H^q(\Sc_\rho\,\B \otimes \O_G(m)) \neq 0$, we have   
$q=l_\gamma$ (in the notation of \cite{ego2}). But now   
\[ \rho_1 \le (d-i+1)(i-r+1) \le d-t+1,\]
so it is immediate that $l_\gamma = 0$ if $m \ge 0$. If $m <0$, then   
$l_\gamma$ can only be a multiple of $d-t+1$. If $t =1$, then   
necessarily $l_\gamma = d+1 = \dim G$. If $t > 1$, then   
$l_\gamma < \dim G$ is possible only if $\rho_1 = d-t+1$, i.e., only   
if C2 is an equality. The lemma is proved, and the proof of the main theorem   
is complete. \qed

\begin{Example} \rm   
\begin{itemize}
\item  
The data $(t,d,i,r)=(3,16,13,11)$ satisfy C1 and strict C2.   
Set-theoretically $\Le(13,11) = \Le(\uh)$ for   
\[  \begin{array}{ccccccccccc}
\uh = 1 & 2 & \dots & 11 & 12 & 13 & 11 & 9 & 6 & 3 & 0.  
\end{array} \]   
Thus $\Le(13,11) \subseteq \P (\bigwedge^3 S_{16})$ is ACM.  
\item  
Similarly the data $(t,d,i,r)=(5,32,28,24)$ satisfy C1 and strict C2.  
In this case $\Le(28,24) = \Le(\uh)$ for
\[  \begin{array}{ccccccccccc}
\uh = 1 & 2 & \dots & 27 & 28 & 24 & 20 & 15 & 10 & 5 & 0.
\end{array} \]  
\end{itemize} \end{Example}

\begin{Example} \rm Choose an integer $s$ such that   
$t \le s < \frac{t(d+1)}{(t+1)}$. Then the data $(t,d,s,s)$ satisfy C1,C2,  
in fact C2 is strict unless $s=t$.  
Set-theoretically $\Le(s,s) = \Le(\uh)$, where $\uh$ is  
defined as in the proof of Theorem \ref{theorem.dimsigma}.  

If $t=1$, then this is the function $\uh_s$ defined in  
\S \ref{goren.2}. By the Gruson-Peskine theorem, we then know that  
$\Le(s,s) = \Le(\uh_s)$ as schemes.  
We do not know if this remains true for $t > 1$.  
\end{Example}   

\begin{Example} \rm   
Let $t=2$ and choose integers $i,d$ such that $i \ge 5$ and  
$3i = 2d + 1$. Then $(2,d,i,i-1)$ satisfy C1,C2, the latter being strict  
iff $i > 5$. (We recover Example \ref{example.lascoux} for $i=5$.)  
In this case, $\Le(i,i-1)$ is reducible with two components of  
dimension $3i-7$ (i.e., codimension $4$) each.  

For instance, let $(i,d) = (9,13)$, then $\Le(9,8) =  
\Le(\uh_1) \cup \Le(\uh_2)$, where  
\[  \begin{array}{cccccccccccc}
\uh_1  = 1 & 2 & \dots & 7 & 8 & 9 & 8 & 7 & 6 & 4 & 2 & 0. \\  
\uh_2  = 1 & 2 & \dots & 7 & 8 & 8 & 8 & 8 & 6 & 4 & 2 & 0.
\end{array} \]  
\end{Example}

\begin{Example} \rm  
The data $(t,d,i,r)=(3,14,11,9)$ satisfy   
C2, but not C1. Indeed,  
$\Le(11,9) = \Le(\uh_1) \cup \Le(\uh_2) \cup \Le(\uh_3)$, where \
\[  \begin{array}{cccccrrccccc}
\uh_1  = 1 & 2 & \dots & 8 & 9 & 10 & 11 & 9 & 7 & 5 & 3 & 0. \\  
\uh_2  = 1 & 2 & \dots & 8 & 9 &  9 & 9  & 9 & 9 & 6 & 3 & 0. \\  
\uh_3  = 1 & 2 & \dots & 8 & 9 & 10 & 10 & 9 & 8 & 6 & 3 & 0.
\end{array} \]  
The components $\Le(\uh_1),\Le(\uh_2)$ have the expected dimension $27$, but  
$\Le(\uh_3)$ is $28$-dimensional.  
\end{Example}

\begin{Remark} \rm  
By Lemma \ref{lemma.c2}, it is easy to produce examples where C2 holds.  
In contrast, C1 is rather restrictive. (Although, for small values of  
$(t,d)$ it is satisfied more often than not.) It would be worthwhile  
to characterise all sequences $\uh$ such that $\Le(\uh)$ is  
ACM (or projectively normal), but it is unlikely that the technique  
used here can be pushed any further.  
\end{Remark}  

\begin{Remark} \rm \label{remark.char}
Some of the results proved here can be extended to char~$>0$ with appropriate  
care. Replacing $S$ by the divided power algebra (see \cite[Appendix A]{IK})
all results until the beginning of \S \ref{geom.description} remain valid in  
arbitrary characteristic. (The reference to partial differentiation should be  
ignored.)  

All results in \S \ref{geom.description} -- \ref{catalecticant} are valid  
for char~$> d$. In section \ref{free.res} we need to assume  
char $=0$, since (inter alia) Lascoux's result and the  
Borel-Weil-Bott theorem fail to hold in positive characteristic.  
\end{Remark}  

\bibliographystyle{plain}   
\bibliography{../BIB/hema1,../BIB/hema2,../BIB/hema3}

\parbox{7cm}{\small  
Jaydeep V. Chipalkatti \\ (jaydeep@math.ubc.ca) \\  
Dept of Mathematics \\ 1984, Mathematics Road \\
University of British Columbia \\  
Vancouver, BC V6T 1Z2 \\
Canada.}
\parbox{7cm}{\small  
Anthony V.~Geramita \\ (tony@mast.queensu.ca, geramita@dima.unige.it) \\  
Dept of Mathematics and Statistics, \\
Jeffery Hall, Queen's University, \\  
Kingston, ON K7L 3N6 \\ Canada. \\  
\hspace*{2cm} and \\  
Dipartimento di Matematica \\  
Universit{\'a} di Genova \\  
Genova, Italia.}  
\end{document}